\documentclass{siamltex}
\usepackage{multirow}
\usepackage{graphicx}
\usepackage{mathrsfs}
\usepackage{color}
\usepackage{algorithm,algorithmic}
\usepackage{amsmath,amsfonts,amssymb,graphicx}%
\usepackage{bbm}
\usepackage{subfigure}
\usepackage{hyperref}
\usepackage{todonotes}

\newtheorem{remark}[theorem]{Remark}
\newtheorem{example}[theorem]{Example}

\def\Range{{\rm Range}}
\def\span{{\rm span}}
\def\orth{{\rm orth}}

\def\red{{\rm red}}

\def\RR{\mathbb{R}}
\def\CC{\mathbb{C}}

\def\itmax{{\tt itmax}}

\def\RK{{\mathcal{RK}}}

\def\IRKA{{\sc irka}}
\def\RIRKA{{\sc r-irka}}
\def\TRIRKA{{\sc tr-irka}}

\def\tolouter{\texttt{tol}_{\rm outer}}
\def\tolinner{\texttt{tol}_{\rm inner}}
\def\tolIRKA{\texttt{tol}_{\rm IRKA}}
\def\linsolvers{\xi_{\rm lin}}
\def\finalsubspace{\ell^{\rm fin}}

\def\chichiR-IRKA{\chi\chi^{\rm rel}_{\rm R-IRKA}}

\def\H2{$\mathcal{H}_2$}
\def\Hinfty{$\mathcal{H}_\infty$}

\title{	A reduced-IRKA method for large-scale \H2-optimal \\model order reduction
 \thanks{Version of \today.}}

\author{Yiding Lin\thanks{School of Mathematics, Southwestern University of Finance and Economics, 555 Liutai Road, Chengdu 611130, China({\tt yiding.lin@gmail.com}).}
        \and Valeria Simoncini \thanks{Dipartimento di Matematica, Universita di Bologna, Piazza di Porta S. Donato, 5, 40127 Bologna, Italy, and
IMATI-CNR, Pavia, Italy.   ({\tt valeria.simoncini@unibo.it}).}}

\begin{document}
\maketitle

\begin{abstract}
The \H2-optimal Model Order Reduction (MOR)  is  one of the most significant frameworks
for reduction methodologies for linear dynamical systems.  In this context, the Iterative Rational
Krylov Algorithm (\IRKA) is a well established method for
computing an optimal projection space of fixed dimension $r$,
when the system has small or medium dimensions.
However, for large problems the performance of \IRKA\ is not satisfactory.
In this paper, we introduce a new rational Krylov subspace projection method with conveniently
selected shifts, that can effectively handle large-scale problems.
The projection subspace is generated sequentially, and
the \IRKA\ procedure is employed on the projected problem to produce 
a new optimal rational space of dimension $r$ for the reduced problem, 
and the associated shifts.
The latter are then injected to expand the projection space. Truncation of older
information of the generated space is performed to limit memory requirements.
Numerical experiments on benchmark  problems illustrate the effectiveness of the new method.

\end{abstract}

\begin{keywords}
\H2-optimal model order reduction, rational Krylov subspace method, \IRKA,
transfer  function
\end{keywords}

\begin{AMS}
34C20, 41A05, 49K15, 49M05, 93A15, 93C05, 93C15
\end{AMS}

\pagestyle{myheadings}
\bibliographystyle{siam}

\markboth{Yiding Lin and Valeria Simoncini}{A reduced-\IRKA\ method for large-scale \H2-optimal  model order reduction} 

\section{Introduction}
We are interested in reduction techniques for the following
classical linear dynamical system
\begin{equation}\label{eqn:mainsys}
	\left\{
	\begin{aligned}
		E\frac{dx(t)}{dt}&=Ax(t)+Bu(t),\\
		y(t)&=C^Hx(t),\\
	\end{aligned}
	\right.
\end{equation}
with $A \in \mathbb{C}^{n \times n}$  $c$-stable,
 and $B\in \mathbb{C}^{n \times m}$, 
 $C\in \mathbb{C}^{n \times p}$, 
aimed at satisfying optimality conditions
of the
associated transfer function  
\begin{equation}\label{eqn:transfer_functions}
h(\mathfrak{s})=C^H(\mathfrak{s}E-A)^{-1}B.
\end{equation}
In the case that $m=1=p$, we will use $b=B$ and $c=C$.

A large number of techniques have been explored to reduce the model dimension 
while maintaining the principal properties of the original system~\cite{MR2155615},%
\cite{doi:10.1137/1.9781611976083},\cite{BCOW.17}. The main idea consists of determining
a convenient subspace (or pairs of subspaces) onto which to project the original system.
The obtained new model has the same structure as in (\ref{eqn:mainsys}), it possibly
retains many of the original dynamics features, but it has significantly
reduced dimensions.  A leading position as effective approximation space
has been taken by rational Krylov subspaces. 
For a chosen dimension $k$ and given parameter set $\mathbb{S}=\{s_1, \ldots, s_k\}$ of
distinct values $s_j$ (shifts),
the rational Krylov subspace associated with $A$, $E$, $B$ is given by
\begin{equation}\label{eqn:pure_RK_definition}
\RK_k(A,E,B,\mathbb{S}):=
range([(A-s_1E)^{-1}B,(A-s_2E)^{-1}B ,\ldots ,(A-s_kE)^{-1}B]) .
\end{equation}
The easiness of the implementation 
together with key interpolatory properties
have led to the widespread use of these projection spaces 
in the design of model order reduction (MOR) strategies.
Indeed, for specifically generated parameters, 
the resulting reduced problem satisfies the Hermitian interpolatory conditions,
which are necessary conditions for the reduced model to be 
\H2-optimal \cite{doi:10.1137/1.9781611976083}.
In the now classical paper \cite{MR2421462}, the authors proposed
 a very successful approch towards this goal, 
called Iterative Rational Krylov Algorithm (\IRKA). This algorithm
determines ideal parameters satisfying these interpolatory
conditions, and since then it
has become the reference method for small and medium size
linear dynamical systems reduction methodologies.
%
\IRKA\ possesses numerous variants\cite{4434939,5400605,6426344,MR3322841,MR3923432,MR2401633,MR4734558} and has been generalized to dealing with nonlinear systems\cite{MR3023455,MR3809536,MR3348125}.  
Further references can be found in the recent works \cite{MR4180031,10740341} and
in the book \cite{doi:10.1137/1.9781611976083}.

Although a considerable number of effective MOR methods have been 
developed, see, e.g., \cite{MR3805855,MR2610803,MR2858340}, 
the \H2-optimal MOR method is distinguished by its optimality properties.
A possible weakness of \IRKA\ is that the number of parameters, which also corresponds
to the reduction space dimension $\ell$, does not generate a nested procedure.
As a consequence, the information obtained in computing 
$(\ell-1)$-order \H2-optimal MOR  is not used for computing $\ell$-order \H2-optimal MOR.
In contrast, the shifts of projection methods, see for instance 
\cite{MR3805855,MR2610803,MR2858340} are nested, thus allowing the space to grow for
better quality. As a result, these procedures
 generally require less CPU time than \IRKA. 
However, the obtained shifts are not \H2-optimal (nor \Hinfty-optimal, in the sense
of \cite[section~5.3]{MR2155615}).

\IRKA\ aims to determine the ideal reduction space by 
generating  a sequence of rational Krylov subspaces of fixed dimension, where the 
set $\mathbb S$ of parameters  used
for their construction is refined at each sequence iteration.
In typical experiments, \IRKA\ behaves like a fixed point iteration 
and converges linearly \cite{MR2924212,10740341}. 
The computational cost may be high because of the repeated
orthogonal bases construction and because of the expensive system solves 
with the matrices $A- s_j E$.
A number of strategies have been taken to develop faster  
algorithms for \H2-optimal MOR\cite{MR4188844,MR3319848,MR4180031}.
As already mentioned, except for \IRKA-based methods,
algorithms usually accumulate the 
obtained shifts \cite{MR2858340,MR4467544,MR2837488,MR3735291}, and
reformulate the problem as that of finding a subset of (quasi) \H2-optimal shifts by
using classical eigenvalue procedures based on the rational Krylov space
 \cite{MR1429696,MR1618804}.
The accumulation of all obtained shifts results in the faster convergence of the algorithm.

In light of these considerations, we derive
 a rational Kyrlov subspace reduction method (in the following \RIRKA) for 
approximating the \H2-optimal shifts.
Our procedure generates a Rational Krylov subspace of dimension $r$, 
then determines \IRKA\ optimal
parameters of the reduced model of size $r$. Then, instead of restarting the process, expands the
rational space with the newly computed parameters as shifts, so as to have a space of dimension $2r$.
Then, a new reduced problem is optimally solved with \IRKA\ by determining
$r$ new parameters, which are then used to further expand the basis, and
so on.
Hence, the reduced problems are classified as 
small size \H2-optimal MOR problems, which can be solved by  \IRKA.
Intuitively, \H2-optimality is ensured in the projected problem; the shifts' inclusion in expanding the
projection space will provide additional information to the original data to further expand the space
towards the sought after ideal reduction space of dimension $r$.

Our methodology can be regarded as a generalization of 
standard RKSM for solving eigenvalue problems, where the eigenfinder for
the reduced problem is replaced by the \H2-optimality finder, \IRKA.
This greedy strategy is similar to other MOR strategies for various
optimization problems associated with dynamical systems, see, e.g.,
\cite{MR3735291,MR4113069,MR4813186,MR4533502,MR4061624,MR4467544,MR3180856}.
Our method differs from previously developed projection schemes in that at
each \RIRKA\ iteration, the next additional block of $r$ rational Krylov vectors, and
not just one, will be
computed using the newly extracted parameters, in a dynamic manner. Moreover,
truncation of the older basis blocks will help control the memory requirements.


Finally, we would like to mention that our approach is based on the original formulation
of \IRKA, which requires the access to a state-space realization of the transfer function
$h(\mathfrak{s})$.
Another viewpoint that has recently been explored consists of using a Loewner-matrix framework,
which only requires the evaluation of the transfer function at 
$\mathfrak{s}\in {\mathbb C}$; for a more detailed discussion
and proper references, we point to the presentation and the
approach in \cite{6426344},\cite{MR3319848}. 
Nonetheless, obtaining a state-space realization
may still valuable, providing additional information towards the stability analysis of the system.

Here is a synopsis of the paper.
In Section \ref{sect:optimality_condition}, we review the  related results about the  \H2-optimal MOR problem.
We  establish  our algorithms and provide implementation details in Section~\ref{sect:projection_algorithms}.
A comparison between \IRKA\ and  \RIRKA\ is provided in Section~\ref{sect:numerical_experiments}.
We draw our final remarks in Section \ref{sect:conclusion}.

\vskip 0.1in
{\bf Notation}:
For $s\in {\mathbb C}$, we denote with $\bar s$ its conjugate.
For $V\in \CC^{m\times n}$, its conjugate transpose is denoted by $V^H$.  
Unless otherwise specified, $\| \cdot\|$ denotes the Euclidean norm for vectors.
The operation $\orth(V)$ produces a matrix whose columns are
an orthonormal basis of $\Range{(V)}$. 
Matlab notation will be employed whenever feasible.
Using standard terminology in control,
Single-Input / Single-Output (SISO) systems are characterized by single columns in $B$ and $C$,
in which case we will use $b\equiv B$ and $c\equiv C$. Multicolumn $B$ and $C$
correspond to Multi-Input / Multi-Output (MIMO) systems.

\section{The \H2-optimal MOR problem}\label{sect:optimality_condition}
To introduce the reduced model, let us first consider a SISO system,
so that $b\equiv B$ and $c\equiv C$.

A reduction technique tries to identify
two full rank matrices $V, W \in \RR^{n\times r}$ that define
the following reduced system 
\begin{equation*}
\left\{
\begin{aligned}
		W^HEV\frac{d}{dt} \tilde x(t)&=W^HAV \tilde x(t)+W^Hbu(t),\\
		y(t)&=c^HV\tilde x(t) ,\\
\end{aligned}
\right. 
\end{equation*}
so that the new reduced system maintains all relevant information
of the
original dynamical system (\ref{eqn:mainsys}). The reduced system is
associated with the transfer function
\begin{equation}\label{eqn:Rtransfer_functions}
\tilde{h}(\mathfrak{s})=c^HV(\mathfrak{s}W^HEV-W^HAV)^{-1}W^Hb,
\end{equation}
whose poles are the eigenvalues of $(W^HV)^{-1}W^HAV$.
The quality of the reduced system can be measured in terms of the
error between the two transfer functions $h$ and $\tilde h$, in some norm.
In this paper we consider the \H2-norm: given the space 
$\mathcal{H}_2 =
\{g : g \mbox{ analytic in } {\mathbb C}_+, 
 \,   \sup\limits_{x>0} \int_{-\infty}^{+\infty}
		|g(x+\iota y)|^2dy<\infty\}$,
equipped with the inner product
$$
\langle g, h\rangle_{\mathcal{H}_2}
:=\frac{1}{2\pi} \int_{-\infty}^{+\infty}
		\overline{g(\iota\omega)}
		h(\iota\omega)d \omega,
$$
we will use the corresponding norm $\|g\|_{\mathcal{H}_2}$.
Hence, the \H2-norm of the model error associated with the reduced
transfer function $\tilde h$ 
defined in (\ref{eqn:Rtransfer_functions})
is computed by means of the
following relative quantity\footnote{The inner product can 
be computed by solving Sylvester 
equations~\cite[Lemma 2.3]{MR2421462} or by evaluating the 
residues of the transfer functions~\cite[Lemma 2.4]{MR2421462}. 
Within Matlab, this computation can be conveniently performed using
 the function \texttt{norm(sys,2)}.}
\begin{eqnarray}\label{eqn:sigma}
\sigma(\tilde h)=
\frac{\|h-\tilde h\|_{\mathcal{H}_2}^2}{\|h\|_{\mathcal{H}_2}^2}.
\end{eqnarray}
The \H2-optimal MOR problem consists of determining $\tilde h$ that solves
the following optimization problem
\begin{eqnarray} \label{eq:H_2_problems}
\|h(\mathfrak{s})-\tilde{h}(\mathfrak{s})\|_{\mathcal{H}_2}&=&
\min_{\substack{ {\rm dim} (\tilde{h}_*(\mathfrak{s}))=r\\
\tilde{h}_*(\mathfrak{s}): ~ {\rm stable}} }
	\|h(\mathfrak{s})-\tilde{h}_*(\mathfrak{s})\|_{\mathcal{H}_2}.
\end{eqnarray}
where the conditions 
``${\rm dim} (\tilde{h}_*(\mathfrak{s}))=r$,
$\tilde{h}_*(\mathfrak{s}): ~ {\rm stable}$'' mean
that the reduced system is stable and has dimension $r$.

A possible solution can be obtained via a
moment matching approach, which consists of
imposing interpolation conditions, and this is described next.
In passing, we observe that when the system does not possess the 
state-space representation, an alternative approach is 
to apply the Transfer Function \IRKA\  \cite{6426344}.

In the following theorem we recall the 
 Meier-Luenberger 
first-order necessary interpolation conditions (\cite{1098680})
for an approximate solution to the \H2-optimal MOR problem 
(\ref{eq:H_2_problems}). 

\vskip 0.1in
\begin{theorem} {\rm \cite[Theorem 3.4]{MR2421462}}
 \label{H_2_optimal_condition}
Suppose $A$ is stable.
Let $\tilde{h}$ be a local minimizer of problem 
\eqref{eq:H_2_problems}.
Suppose  that $\tilde{h}$ has simple poles at $\lambda_i$, $i=1,2,\dots,r$.
Then, it holds that
$h(-\overline{\lambda_i})=\tilde{h}(-\overline{\lambda_i})$, 
$h'(-\overline{\lambda_i})=\tilde{h}'(-\overline{\lambda_i})$ for $i=1,2,\ldots r$.
\end{theorem}
\vskip 0.1in

The classical \IRKA\ algorithm determines $V, W$ spanning rational Krylov
subspaces 
$\Range (V)=\RK_r(A,E,b,\mathbb{S})$ and 
$\Range (W)=\RK_r(A^H,E^H,c,{\mathbb{S}})$
such that the 
corresponding poles aim to satisfy the interpolation conditions of Theorem ~\ref{H_2_optimal_condition}
\cite{MR2421462}.
We explicitly observe that the shifts are the interpolation nodes 
of $h(\mathfrak{s})$, thus
the conditions for the first-order derivative require the 
shifts for the left and right subspaces to be the same. 
At convergence it holds that $s_i=-\lambda_i,  i=1,2,\ldots r$,
where the $\lambda_i$ are the eigenvalues of the reduced system matrix
$(W^H EV )^{-1} W^H AV$ \cite{MR2421462}.
If all $A$, $b$ and $c$ are real, then the shifts set should be
closed with respect to conjugation, so that only real or 
conjugate pairs shifts are involved.  
The typical \IRKA\ iteration for the SISO case is sketched as follows,

\vskip 0.1in
\hskip 0.3in Select set of shifts ${\mathbb S}_1$

\hskip 0.3in For $k=1, \ldots,$ 

\hskip 0.4in Construct $V, W$ spanning 
$\RK_r(A,E,b,\mathbb{S}_k)$, $\RK_r(A^H,E^H,c,\overline{\mathbb S}_k)$, resp

\hskip 0.4in Compute eigenvalues $\{\lambda_i\}_{i=1\ldots, r}$  
of $(W^H EV )^{-1} W^H AV$

\hskip 0.4in Set $\mathbb{S}_{k+1}=\{-\lambda_1, \ldots, -\lambda_r\}$ 

\hskip 0.4in If satisfied then stop
\vskip 0.1in

The complete algorithm is recalled in Algorithm~\ref{alg:IRKA} in Appendix 1.
The behavior of \IRKA\ has been thoroughly 
researched \cite{10740341,MR4180031,MR2924212}.

\section{The reduced \IRKA\ method} \label{sect:projection_algorithms}
Although very effective for small to medium size problems, \IRKA\ becomes
expensive for large scale matrices, where the main cost is the solution
of very many linear systems to construct the two Krylov subspaces at each
\IRKA\ iteration.
A key fact in this respect is that,
to maintain the dimensions of the rational Krylov subspace always 
at most equal to $r$, \IRKA\ discards the bases computed in previous iterations.
On the other hand, it has been noticed that rational Krylov subspaces
adapt the shifts as the iterations proceed \cite{MR2858340},
allowing one to determine a good approximation to the transfer function.
Hence, our idea is to allow the rational Krylov subspaces to grow, and
seek an optimal set of $r$ shifts by means of a  ``reduced'' \IRKA\ 
applied to the (small) projected matrices in the generated subspaces.
The procedure comprises the following steps, where \IRKA($r$) means that \IRKA\
generates two bases of dimension $r$, together
with the corresponding approximate optimal $r$ shifts.

\vskip 0.1in

\hskip 0.2in Select initial $\mathbb{S}_0$

\hskip 0.2in Compute $\widehat V, \widehat W$
generating $\RK_{2r}(A,E,b,\mathbb{S}_0)$, 
$\RK_{2r}(A^H,E^H,c,\overline{\mathbb S}_0)$, resp

\hskip 0.2in For $j=1, \ldots,$ 

\hskip 0.3in Project $A, E, b, c$: 
 ${A}_{\rm red}=\widehat W^HA \widehat V$,  
$E_{\red}= \widehat W^H E\widehat V$,  
${b}_{\red}=\widehat W^Hb$, ${c}_{\red}=\widehat V^Hc$

\hskip 0.3in  Determine $\mathbb{S}_j$ by applying \IRKA($r$) to 
$({A}_\red, {E}_\red, {b}_\red, {c}_\red)$

\hskip 0.3in  If satisfied then stop

\hskip 0.3in  Expand $\widehat V, \widehat W$ by computing
$\RK_r(A,E,b,\mathbb{S}_j)$, $\RK_r(A^H,E^H,c,{\overline{\mathbb S}}_j)$, resp

\vskip 0.1in
We explicitly observe that
the projection subspaces Range$(\widehat V)$ and Range$(\widehat W)$ are constructed 
by accumulating all computed bases during the iterations. In practice, the dimension of
these spaces will likely be less than $r$ times the number of iterations, hence the
orthogonalization of the extended basis needs be carried out at the end of each iteration.
The major advantage of this procedure is that \IRKA\ is now applied to 
significantly smaller matrices, allowing one to limit the computational costs.
On the other hand, we already remarked that 
the dimensions of
the two matrices $\widehat V, \widehat W$ grow by $r$ vectors (at most), so 
that the dimensions of the reduced problem grow correspondingly. 
In section~\ref{sect:R-IRKAD} we describe how to alleviate this problem, in case memory 
allocation constraints arise.


\begin{algorithm}[tbh]
\caption{A Reduced-\IRKA\ (\RIRKA) method for the \H2-optimal MOR problem.
\label{alg:R-IRKA_H2_optimal}}

\begin{algorithmic}[1]
\REQUIRE $A,E \in \mathbb{C}^{n \times n},b,c\in \mathbb{C}^n$,  $r$ (\IRKA\ dim.), 
Option ($\#$ for initialization)
\STATE Get initial bases:	
\begin{enumerate}
\item[\hskip 0.3in{Option 1:}] 
{Approx orth eigenbasis $\widehat V$
 of $2r$ smallest eig's of $(A,E)$},
 $\widehat W = \widehat V$;	
\item[{Option 2:}] 
Obtain $\widehat W$  and  $\widehat V$ from the   known 
 approximate system (*);
\item[{Option 3:}] 
$\widehat V = \orth(\texttt{randn}(n,2r)), \widehat W = \orth(\texttt{randn}(n,2r))$;
\end{enumerate}
\STATE
$\mathbb{S}_0=\texttt{ones}(r,1),$
\FOR{$j=1,2, \ldots, k^{\max} $}
\STATE ${A}_{\rm red}=\widehat W^HA \widehat V$,  
$E_{\red}= \widehat W^H E\widehat V$,  ${b}_{\red}=\widehat W^Hb$, ${c}_{\red}=\widehat V^Hc.$
\STATE Apply \IRKA$(r)$  to find $r$ 
 shifts $\mathbb{S}_k$ for  
${h}_\red(\mathfrak{s})= c_\red ^H(\mathfrak{s}  E_\red - A_\red )^{-1} b_\red $;
\STATE Compute $V$ and $W$ s.t. 
\vskip -.2in
$$
\Range(V)=\RK({A},{E}, b,\mathbb{S}_k,r), \,\,
\Range(W)=\RK({A}^H,{E}^H, c,{\overline{\mathbb{S}}_k},r);
$$
\vskip -.1in
\STATE 	$\widehat V=\orth([\widehat V,V])$, $\widehat W=\orth([\widehat W, W]);$
	\hfill	\%	Expand bases
\hfill		\%	Accumulate  shifts
			\STATE  {\bf if}
$\|\mathbb{S}_k-\mathbb{S}_{k-1}\|_2 /\|\mathbb{S}_k\|_2
			< \texttt{tol}_1$  {\bf then} {\bf break};
\hfill \% Test shift variation
\ENDFOR
	\STATE 		$V=\orth(V),W=\orth(W)$
			
			\ENSURE   \H2-optimal MOR  shifts and subspaces: 
$\mathbb{S}_k$, $\Range({V})$, $\Range({W})$.
\end{algorithmic}
{\footnotesize (*)We will not explore this option in our experiments, although it
may provide an efficient starting point, whenever available. It is also included in \IRKA. } 
\end{algorithm}

A more detailed description of the new method is reported 
in  Algorithm \ref{alg:R-IRKA_H2_optimal}.

A few additional comments on the procedure in Algorithm~\ref{alg:R-IRKA_H2_optimal}
are in order.

The initialization step is similar to the one available for \IRKA\ (Algorithm \ref{alg:IRKA}), in particular for Option 1. Options 2 and 3 are included in case Option 1
is too expensive, depending on the problem dimensions.
At iteration $j=1$, the \IRKA\ procedure is applied to the
initial reduced problem corresponding to Option 1; at later iterations,
that is for $j>1$, the \IRKA\ procedure is applied to the
reduced problem using Option 2, with the available
set of shifts as starting guess.

The inner stopping criterion for \IRKA\ (Algorithm \ref{alg:IRKA}, line 8)
is related to the outer stopping criterion
at line 9, and we required that  $\tolinner \leq \tolouter$.
The accuracy obtained for ${\mathbb S}_k$ in \IRKA\ is indeed crucial for the
overall performance of the method. Both stopping criteria are based on
the discrepancy between successive approximations, which is a common criterion
in the literature.
An alternative would be to evaluate
the discrepancy successive transfer functions in the ${\mathcal H}_2$-norm \cite{MR4188844}.

As $\mathbb{S}_j$ is constrained to real numbers or conjugate pairs, 
the two bases $V$ and $W$ can be computed so as to both contain only real values,
so that real arithmetic can be used throughout.

In terms of computational costs, the most expensive step is again the solution
of the linear systems in forming the two rational Krylov subspaces. Currently
Matlab backslash (\cite{matlab.7})
is used for this step, but iterative methods can be considered as well.
Analogously, the orthogonalization step is performed via the QR function in Matlab.

At convergence, the algorithm yields the components $V, W, \mathbb{S}_k$ of
the  transfer function $h_\red^{{\mbox{\IRKA}}}(\mathfrak{s}):=
c^HV(\mathfrak{s}W^HEV-W^HAV)^{-1}W^Hb$, which
 is the  \H2-optimal reduced model for  the original transfer function $h$ at the
last computed set of shifts $\mathbb{S}_k$. 

We conclude with a remark that motivated our algorithmic development.
A more indepth analysis of convergence will be the topic of future research.

\begin{remark}
Let
 the subspaces ${\mathcal RK}(A, E, b, {\mathbb T}, \ell)$,
 ${\mathcal RK}(A^H, E^H, c, \overline{\mathbb T}, \ell)$
be given, where ${\mathbb T}$ collects all computed sets ${\mathbb S}_k$. Let
$\widehat V, \widehat W$ be the corresponding bases. Assume that the spaces
are optimal, that is, the transfer function
$\widehat h(s)=c^H\widehat W^H (A-s E)^{-1}\widehat V b$
interpolates $h$ at the nodes $t_j \in {\mathbb T}$.
If the optimal R-IRKA $r$ shifts ${\mathbb S}$ generated
using $\widehat V, \widehat W$ are a subset of ${\mathbb T}$, then
the \RIRKA\ transfer function also interpolates the original transfer
function, at the nodes in ${\mathbb S}$.

This intuitive fact comes from the property that polynomials of degree $\ell$
are exactly represented in (rational) Krylov subspaces of dimension
larger than $\ell$.
\end{remark}

\subsection{MIMO version of \RIRKA}\label{sect_MIMO_new_algorithm}
The \IRKA\ algorithm can be generalized to the case of multiple inputs and outputs,
that is to MIMO systems, with $B, C$ having $m>1$ columns each. 
A well studied strategy avoids the use of all inputs and output
columns by using {\it tangential} interpolation, which was introduced in \cite{MR2124150}.
The procedure relies on the prior computation or availability of two 
sets of vectors, ${\mathbbm b}_j, {\mathbbm c}_j$, $j=1,\ldots, r$, which are
used to build the first two rational Krylov subspaces,
\begin{eqnarray*}
&&\span\{(A-s_1E)^{-1}B{\mathbbm b}_1,(A-s_2E)^{-1}B{\mathbbm b}_2 ,\ldots ,(A-s_kE)^{-1}B{\mathbbm b}_k\} ,
\\
&&\span\{(A^H-\bar s_1E^H)^{-1}C{\mathbbm c}_1,(A^H-\bar s_2E^H)^{-1}C{\mathbbm c}_2 ,\ldots ,(A^H-\bar s_kE^H)^{-1}C{\mathbbm c}_k\} .
\end{eqnarray*}
After one iteration, the reduced transfer function is written in terms of a pole/residue
expansion, that is 
$$
h_\red(s) = \sum_{i=1}^r \frac{\widehat {\mathbbm b}_i \widehat {\mathbbm c}_i^H}{s-\lambda_i}.
$$
Setting 
${\mathbbm b}_j = \widehat {\mathbbm b}_j$, 
${\mathbbm c}_j = \widehat {\mathbbm c}_j$ , $j=1,\ldots, r$ gives the tangent combination for
the next iteration. We refer to Chapter 5 in \cite{doi:10.1137/1.9781611976083} for a 
collection of recent related results, together with a rich
analysis of this tangential setting, where interpolation optimality properties are also discussed.
The tangential \IRKA\ method is summarized in
\cite[Algorithm 5.2.1]{doi:10.1137/1.9781611976083}.

Similarly, we can extend Algorithm~\ref{alg:R-IRKA_H2_optimal} to dealing with MIMO systems.
Like for \IRKA, the principal modification in implementing the tangential variant
emerges in constructing the basis matrices of the projection subspaces.
%
%
%
At each iteration
the reduced problem is still defined by the projected matrices
 ${A}_{\rm red}=\widehat W^HA \widehat V$,  
$E_{\red}= \widehat W^H E\widehat V$, and
${B}_{\red}=\widehat W^H B$, ${C}_{\red}=\widehat V^H C$, with
The inner iteration is now the tangential
MIMO version of \IRKA, and all quantities are updated 
accordingly, including the tangential directions.
%
Despite the fact that the whole algorithm also involves the convergence of the tangent directions, only the information of the shifts $\mathbb S_k$ is used for the purpose of stopping the outer iteration.

\subsection{Truncated version of \RIRKA}\label{sect:R-IRKAD}
As the outer iteration proceeds, the dimensions of the projection  subspaces 
undergo a rapid increase due to the accumulation of all the obtained shifts. 
A reduction in memory requirements can be achieved through the 
``deflation'' of some bases vectors.
In fact, we propose to {\it truncate} the basis, by 
retaining only the bases constructed in the last
$\tau$ iterations.
The value of $\tau$ has been determined empirically by a considerable amount of 
experimental testing.
Setting $\tau = 1$ indicates that the obtained shifts are not accumulated, so that
 \RIRKA\ reduces to \IRKA.
While for $\tau=2$ we observed a systematically larger number of outer iterations to reach
convergence,
the choices $\tau = 3,4$ exhibited minimal variations with respect to the non-truncated case,
so that the value $\tau=3$ has been used in our experiments.

Truncation ensures 
 that only the last $3r$ computed vectors for $\widehat V$ and
 $\widehat W$ are retained, where for this discussion we have 
assumed full rank in the basis expansion,
and orthogonalization is performed keeping track of the column ordering.
The term {\it truncated} is reminiscent of a similar strategy used for projection-based
iterative linear system solvers; see, e.g., \cite[sections 6.4.2-6.5.6]{MR1990645}. However, as opposed to the
latter setting, the dimension of the reduced problem is truncated as well.


The truncation process occurs at Line 7 of Algorithm~\ref{alg:R-IRKA_H2_optimal}, which
is replaced by the following step, assuming that $\widehat V$ and $\widehat W$ have $r\cdot k$ columns each:

\vskip 0.1in
\hskip 0.1in Line 7': 
$\widehat V=\orth([\widehat V(:,r (k-2)+1: r k), V])$

\hskip 0.6in $\widehat W=\orth[\widehat W(:,r (k-2)+1: r k), W])$.
\vskip 0.1in

This strategy is preferable to purging columns {\it after} the orthogonalization
in the original Line 7. Indeed, with the chosen strategy we are able to retain
all information contained in the latest computed \IRKA\ basis.
We experimentally observed that this strategy leads to faster convergence than
the truncation after the full orthogonalization, requiring a couple fewer outer
iterations to reach convergence.

\section{Numerical experiments}\label{sect:numerical_experiments}
In this section we report some of our extensive numerical experiments
with the new method \RIRKA\ (Algorithm~\ref{alg:R-IRKA_H2_optimal})
and its truncated variant, which we call \RIRKA($\tau$), described
 in Section~\ref{sect:R-IRKAD}. Our reference algorithm is \IRKA\ (Algorithm~\ref{alg:IRKA}),
although \IRKA\ is designed for medium to large matrices. We are not aware
of other strategies for large scale problems that solve the same optimization
problem.
In the case of a  MIMO problem, the tangential MIMO version of \RIRKA\ 
described in Section~\ref{sect_MIMO_new_algorithm} is used.
All algorithms for the following examples terminate when the stopping criterion
(Line 9 in Algorithm~\ref{alg:R-IRKA_H2_optimal}) 
\begin{equation}\label{eqn:stopcrit}
\chi_k :=\|\mathbb{S}_k-\mathbb{S}_{k-1}\|_2 /\|\mathbb{S}_k\|_2
			< \texttt{tol}_{outer}  .
\end{equation}
The maximal iteration numbers of the outer and inner iterations 
in \RIRKA\ are  ${\itmax}_{\rm outer}=30$
and  $\itmax_{\rm inner}=300$, respectively.
The maximal iteration number of \IRKA\
is ${\itmax}_{\rm IRKA}=300$.
The meanings of the remaining symbols are displayed in Table~\ref{tbl:exlaining_symbols}.

{\it Experimental environment.}
All experiments were carried out in Matlab2021a on a 64-bit notebook 
computer with an Intel CPU i9-11900H and 32GB  memory.
Any data involving random numbers is fixed  
by setting \texttt{rand(`state',0)}, \texttt{randn(`state',0)}
or \texttt{randn(`state',10)}.
The Matlab eigenvalue function
\texttt{eigs} uses a random vector as the initial vector,
hence we set \texttt{rand(`state',0)} before calling this function.

The Matlab code implementing our new algorithm will be made available
by the authors shortly.


\begin{table}[H]
\centering
\caption{   Notation \label{tbl:exlaining_symbols}}
\begin{tabular}{|c|l|}
				\hline
Symbols &\qquad  \qquad \qquad	Explanations	\\
\hline	
SYM & `Yes' for both $E$ and $A$ symmetric\\
\hline
{\#its}&  final  number of iters at termination \\
\hline
\multirow{2}{*}{$\linsolvers$}
& number of linear eqn solves: $\xi_{\rm lin}= 2r \times (\#its) $ \\
&{\footnotesize(For \RIRKA, only outer linear solves are considered)}\\
\hline
\multirow{3}{*}{$\finalsubspace$} & Subspace dimension at termination \\
		&(the total memory allocation is $2n\finalsubspace$: \\
	& $\ell^{\rm fin}_{\mbox{\IRKA}}=r $, \quad $\finalsubspace_{\mbox{\RIRKA}}=2r+r \times (\#its)$, \quad $\finalsubspace_{\mbox{\TRIRKA}}=3r)$ \\
\hline
\end{tabular}
\end{table}

\begin{figure}[htb]
\centering
\caption{Example~\ref{ex:iss} (ISS), $b=B(:,1),c=C(:,2)$. Option 1 is used as starting guess for all methods.
Left: Update in the shift set, $\chi_k$ in 
(\ref{eqn:stopcrit}). Right: Shifts at termination.   \label{figure:SISO_iss_different_order} }
\includegraphics[width=2.2in]{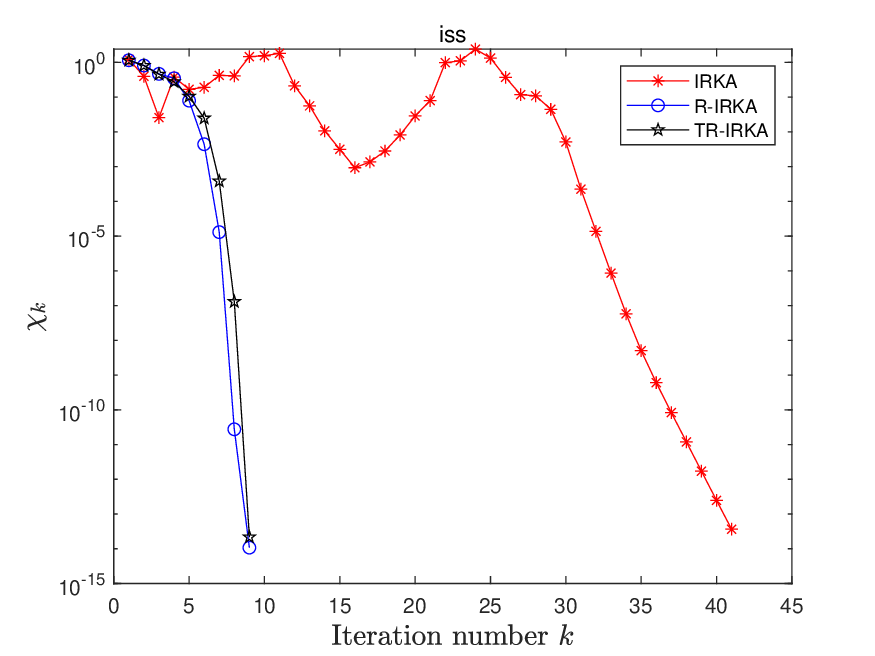}
\includegraphics[width=2.2in]{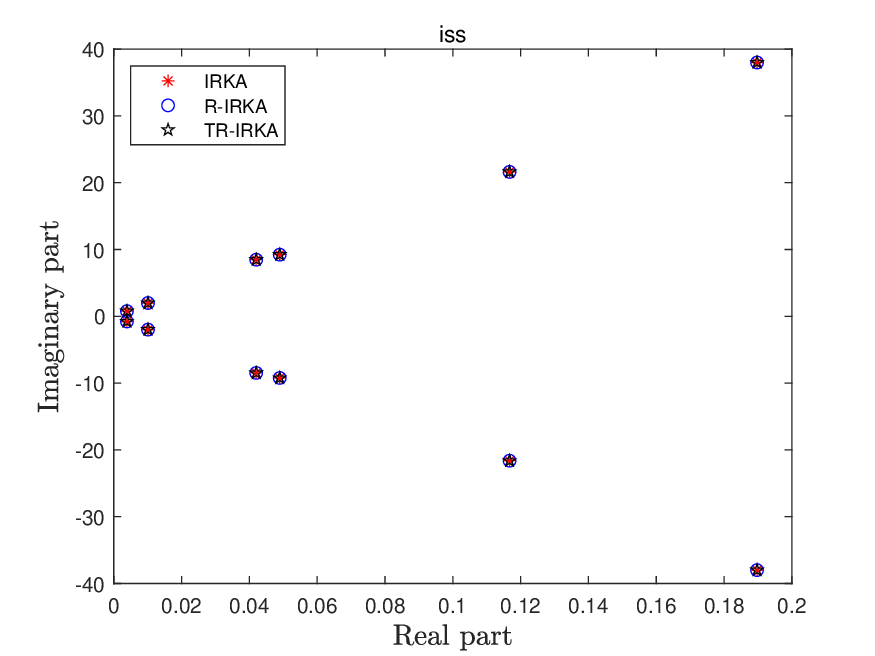}
\\
\center{ \texttt{ $r=12$}, \qquad 
		$\sigma_{\mbox{\IRKA}}=\sigma_{\mbox{\RIRKA}}=
		\sigma_{\mbox{\TRIRKA}}=0.82748$}\\
		\includegraphics[width=2.2in]{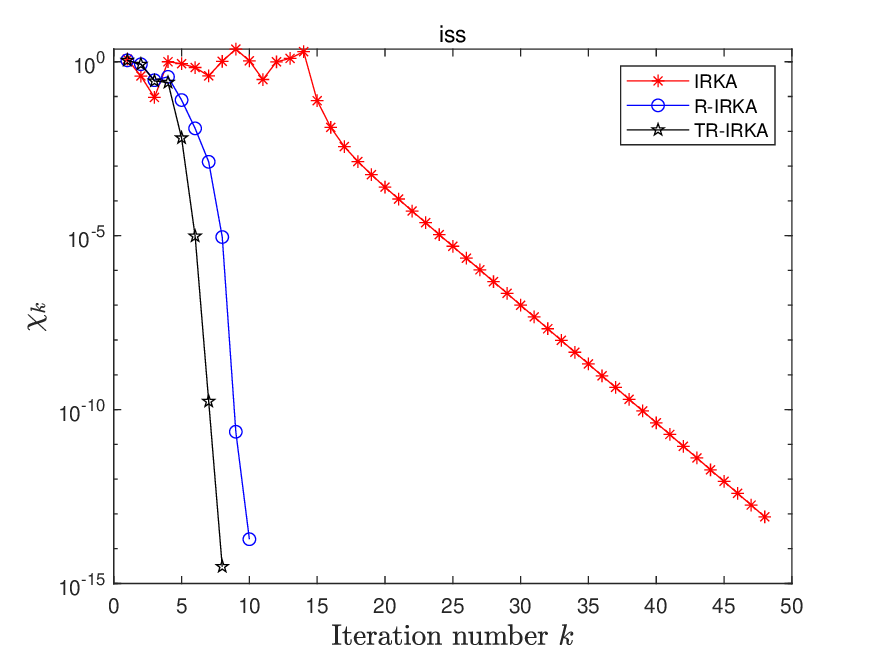}
		\includegraphics[width=2.2in]{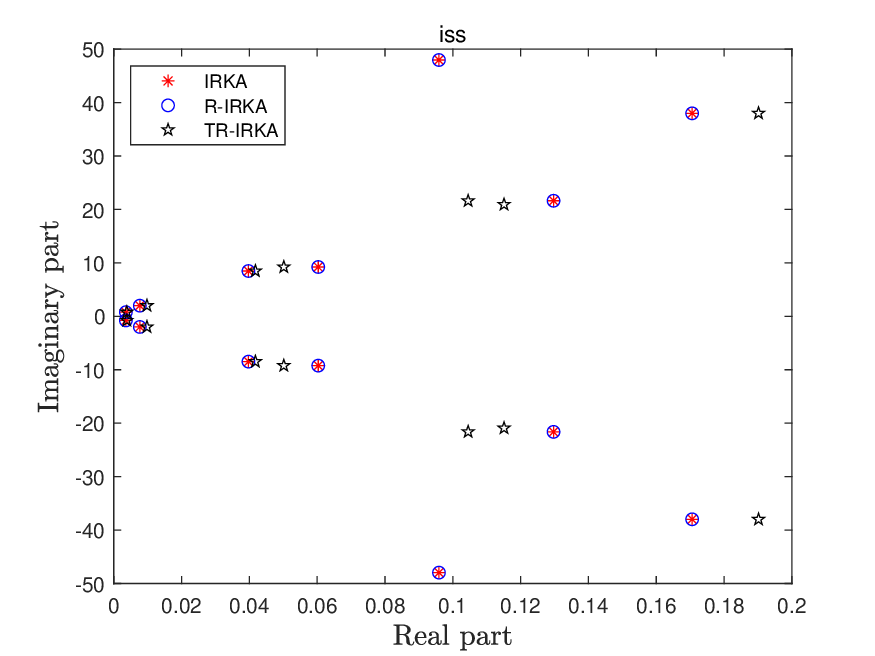}
\\
	\center{\texttt{ $r=14$},  \qquad 
$\sigma_{\mbox{\IRKA}}=\sigma_{\mbox{\RIRKA}}=0.099631,  \sigma_{\mbox{\TRIRKA}}=0.82588$}
\\
		\includegraphics[width=2.2in]{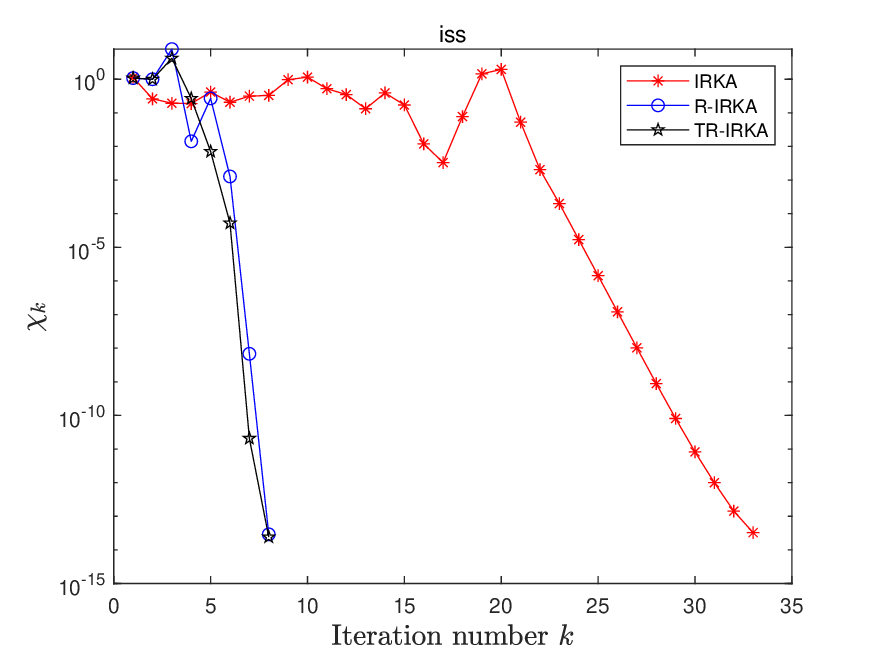}
		\includegraphics[width=2.2in]{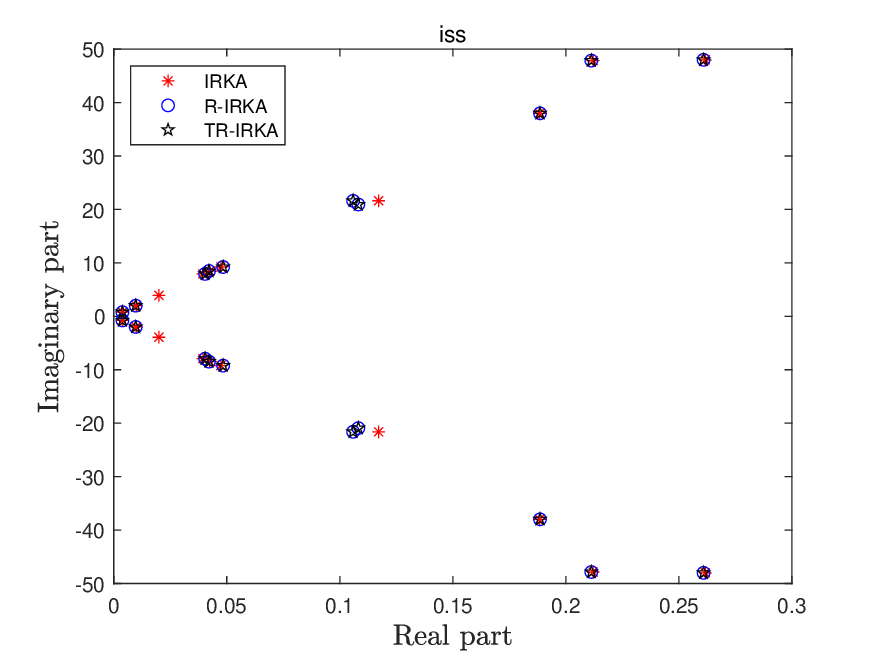}
\\
	\center{ \texttt{ $r=20$}, \qquad 	$\sigma_{\mbox{\IRKA}}=0.0025156, 
\sigma_{\mbox{\RIRKA}}=\sigma_{\mbox{\TRIRKA}}=0.0010077$}
\\	
\end{figure}

\begin{figure}[htb]
\centering
	\caption{Example~\ref{ex:iss}. Option 1 was used for the starting guess of all methods.
Left: Plot of the transfer function and of its
approximations. Right: $|h(\mathfrak{s})-\tilde h(\mathfrak{s})|$.
		\label{figure:SISO_iss_hs_hs_error} }
	
		\includegraphics[width=2.2in]{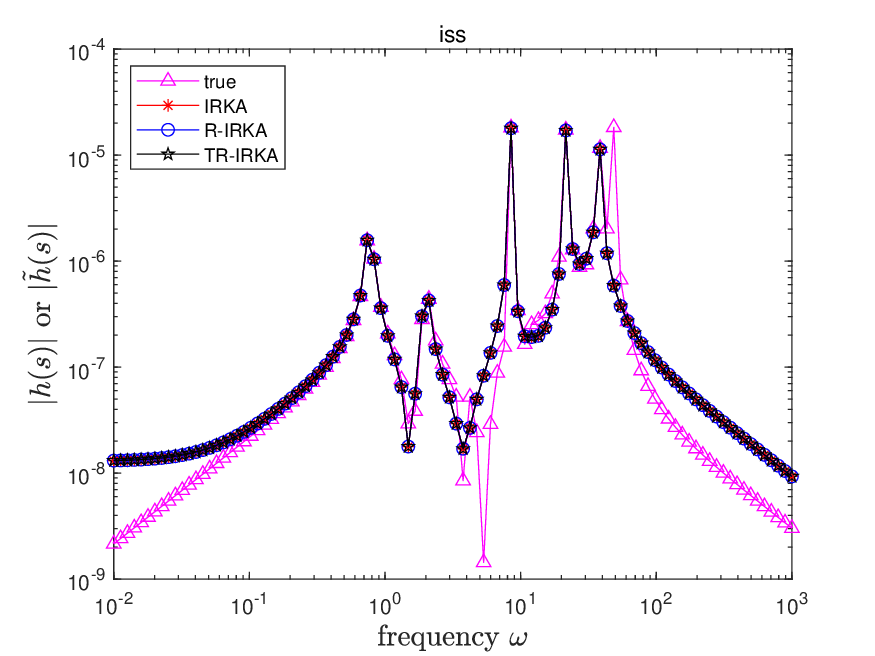}
		\includegraphics[width=2.2in]{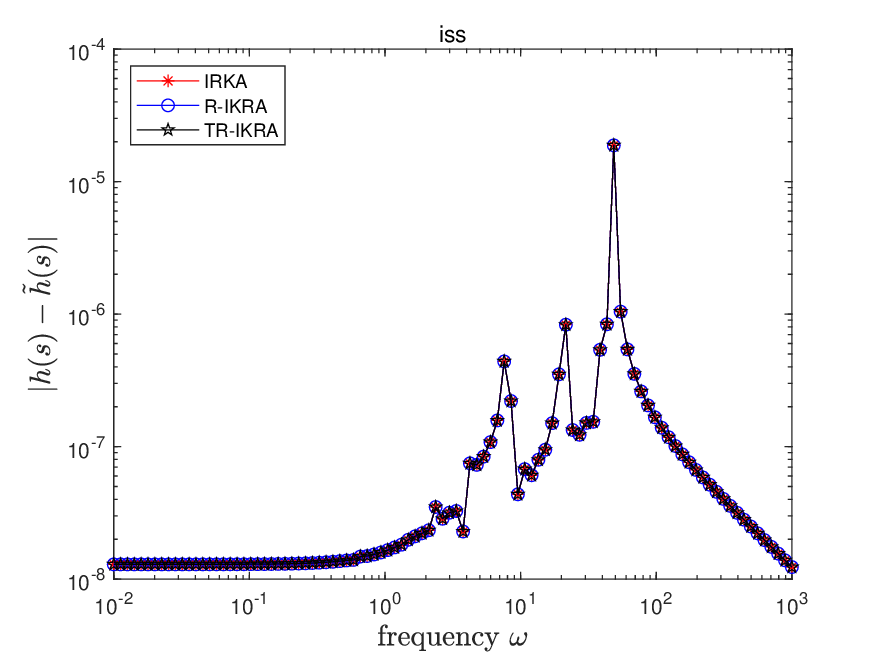}
\\
	\center{ \texttt{ $r=12$}}.
\\
		\includegraphics[width=2.2in]{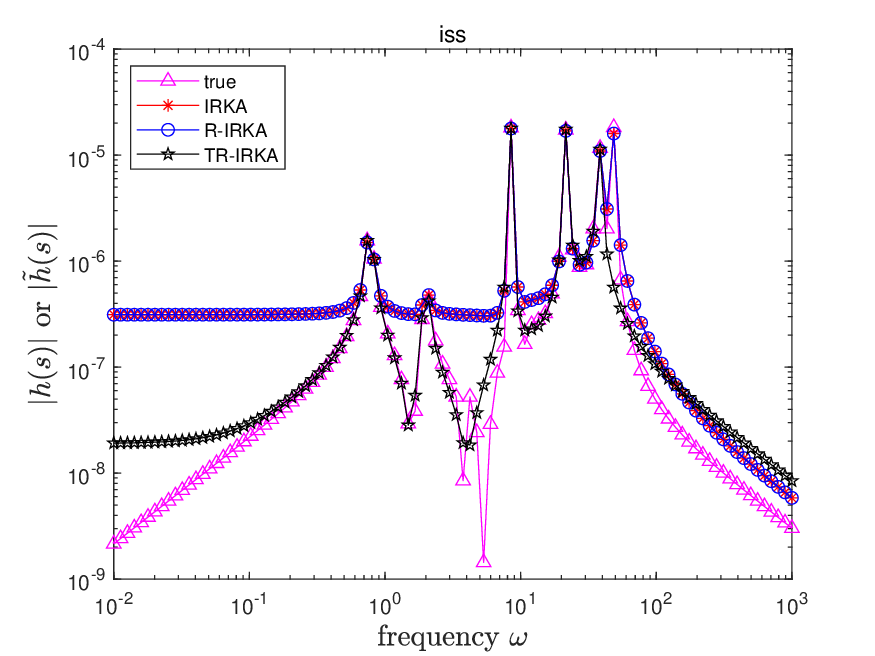}
		\includegraphics[width=2.2in]{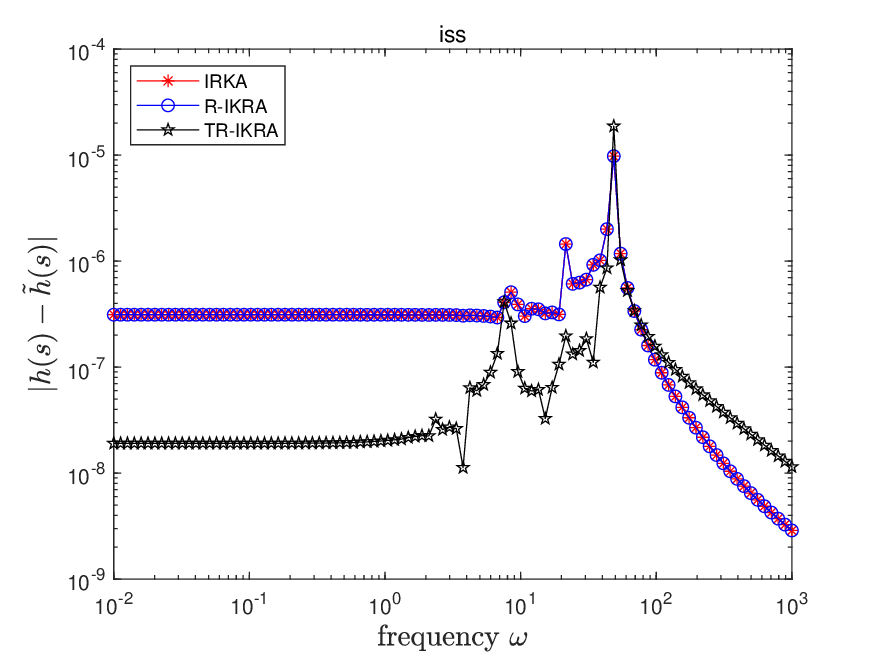}
\\
	\center{	\texttt{ $r=14$}.}
\\
		\includegraphics[width=2.2in]{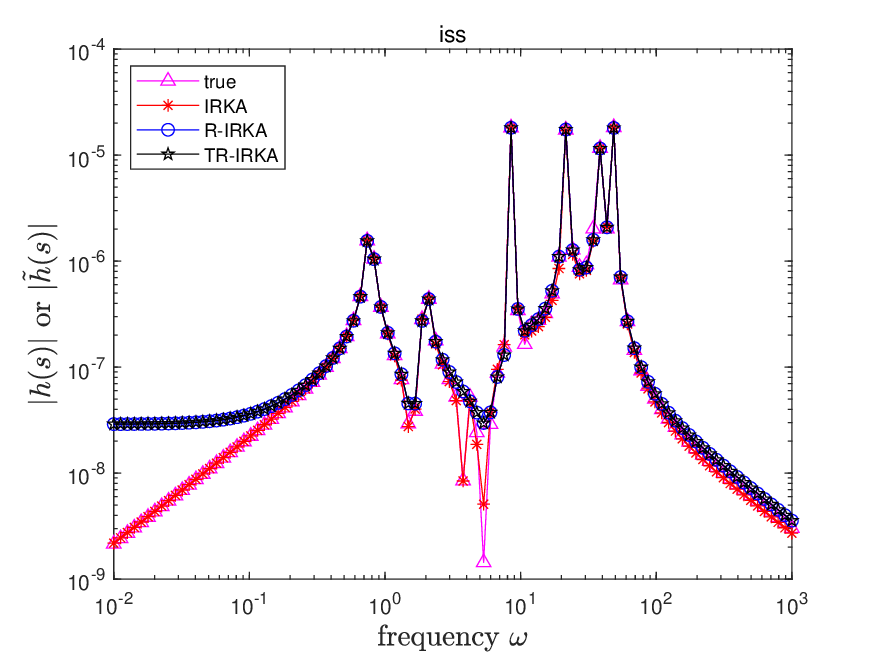}
		\includegraphics[width=2.2in]{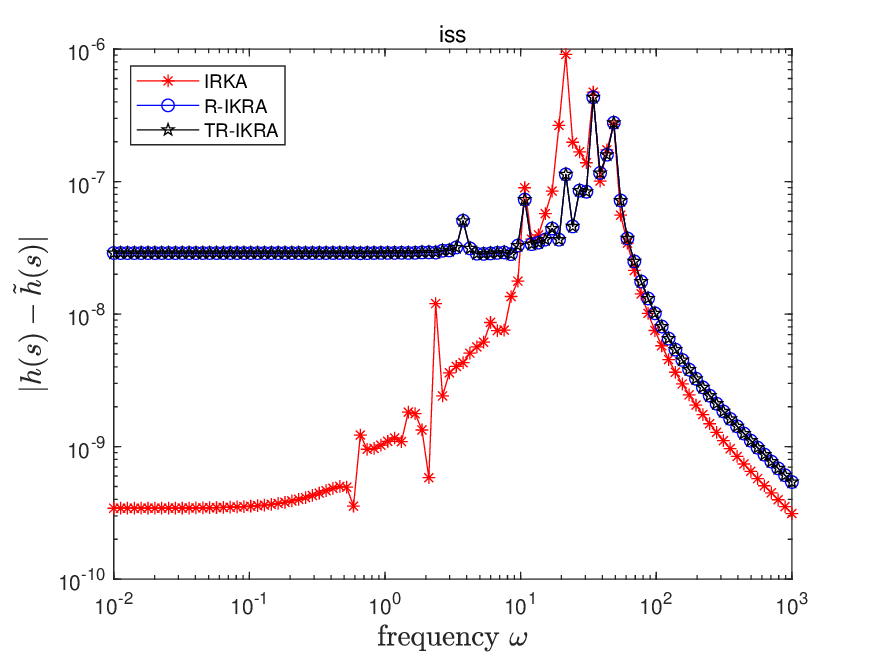}
\\
	\center{\texttt{ $r=20$}.\\}
\end{figure}

\vskip 0.1in
\begin{example}\label{ex:iss}
{\rm
We first consider a SISO small example, the International Space Station (ISS)
benchmark dataset (\cite{Collection2003}), with $A$ of dimension 
270 and both $B$ and $C$ having $3$ columns, and we selected
$b=B(:,1),c=C(:,2)$. 
The problem size allows us to perform a detail error analysis.
In Figure~\ref{figure:SISO_iss_different_order} the performance of the methods is
reported, for $r=12, 14, 20$, with tolerance values
 $\tolinner=5\cdot 10^{-14}$, $\tolouter= \tolIRKA=1\cdot 10^{-13}$.
 The left plots show the quantity $\chi_j$ defined
in the stopping criterion (\ref{eqn:stopcrit}) as the outer iterations proceed.
The right plots display the final shifts obtained by each method.
The ${\mathcal H}_2$-norm (\ref{eqn:sigma}) for all methods and all values of $r$ is also reported.
We note that at termination, the computed shifts satisfy the Meier-Luenberger 
conditions for all methods. We also recall that this condition does not ensure
that global optimizers have been found.

The results in the left plots
show that \RIRKA\ is able to profit from the inner optimal shift selection,
so that the number of new rational Krylov subspace vectors with the original
data is drastically lower than for \IRKA\ (about $2r\cdot 10$ for \RIRKA\ versus 
at least $2r\cdot 35$ for \IRKA).
The right plots show that the for $r=12, 14$ \IRKA\ and \RIRKA\ obtain the same
final shifts. In the truncated version most of the final shifts differed for $r=14$.  
For $r=20$, instead, a few shifts computed with \IRKA\ did not perfectly match
those of the other methods. Note that in this last case,
the ${\mathcal H}_2$-norm is smaller for \RIRKA\ and its truncated variant than with \IRKA. 

Figure~\ref{figure:SISO_iss_hs_hs_error} reports the original transfer function 
and its approximations (left) and the errors
 $|h(\mathfrak{s})-\tilde h(\mathfrak{s})|$ (right).
For all values of $r$, the new method and its truncated variant perform comparably well
with \IRKA\, except possibly for $r=20$ for which better accuracy in some $\omega$ intervals
can be appreciated.


\flushright $\diamond$
}
\end{example}
\vskip 0.1in

\begin{figure}[htb]
	\caption{Example~\ref{ex:largescale}. Data for	\textsc{rail20209}, MIMO case, and \texttt{ $r=11$}
		\label{figure_MIMO_rail20209} }
\centering
		\includegraphics[width=2.3in]{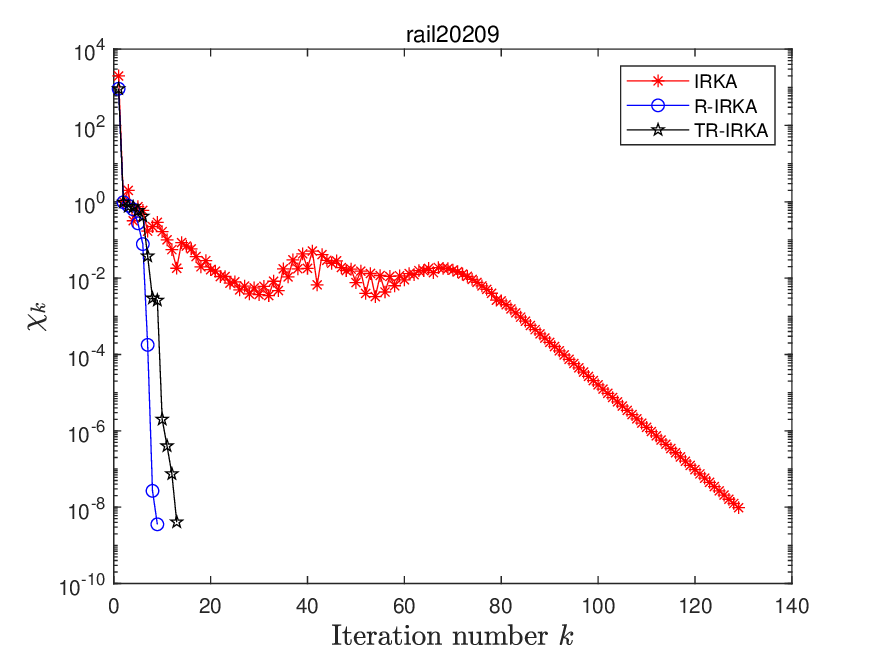} 
		\includegraphics[width=2.3in]{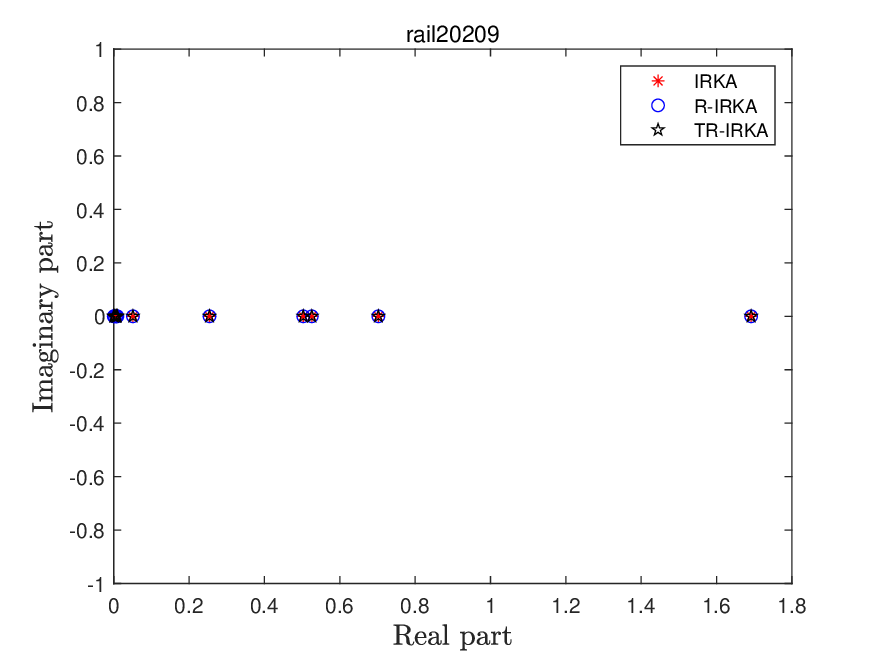}
\end{figure}

\begin{figure}[htb]
\caption{ Example~\ref{ex:largescale}.	Data for 
\textsc{rail20209}, MIMO system, \texttt{ $r=11$}.
Number of inner iterations for \RIRKA\ (left) and for \TRIRKA\ (right) with corresponding CPU time (in seconds).
\label{figure_MIMO_rail20209_T} }
\centering
	
		\includegraphics[width=2.3in]{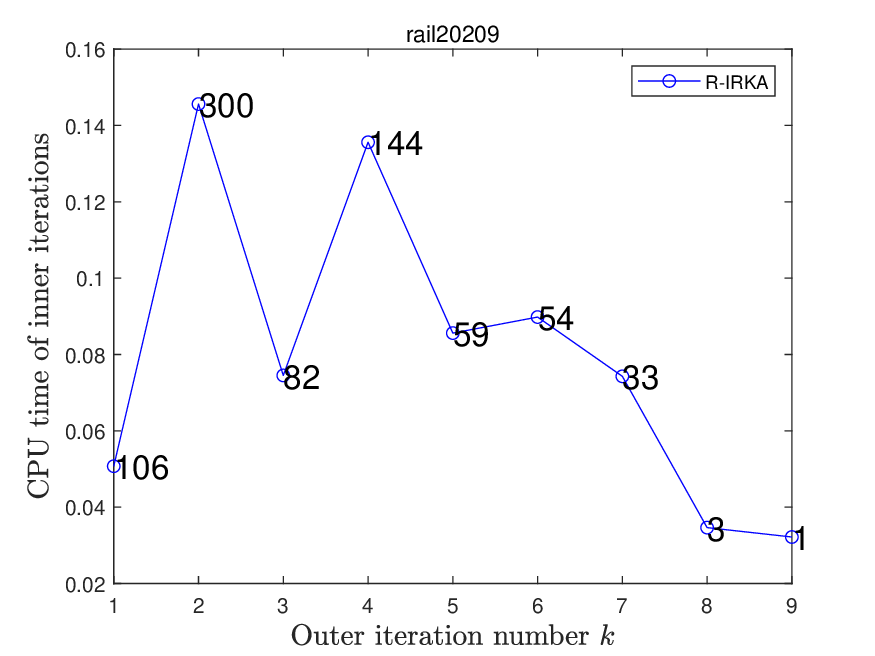} 
		\includegraphics[width=2.3in]{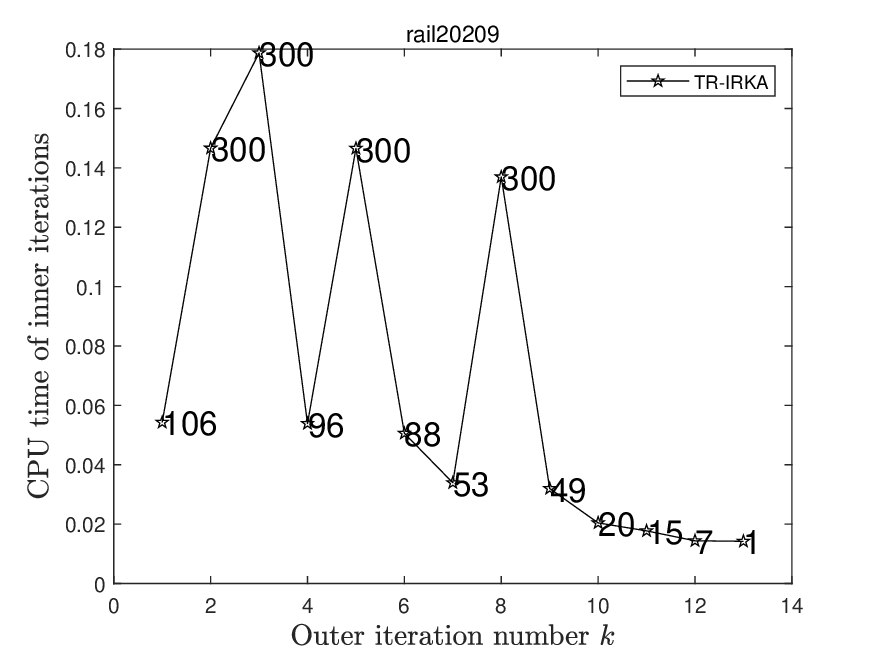} 
	
\end{figure}

\vskip 0.1in
\begin{example}\label{ex:largescale}
{\rm
We consider systems where $A$ is the finite diffence discretization of a two-dimensional
elliptic operator, with zero boundary conditions. The considered operators are
summarized in Table~\ref{tab:matrices_L}.
We set $E=I$ and 
\texttt{randn(`state',0);} \texttt{B=randn(n,2);} \texttt{randn(`state',10);} \texttt{C=randn(n,2)},
specifically, to also  work with MIMO systems.
The remaining problems are obtained from the Oberwolfach collection
\cite{10.1007/3-540-27909-1_11}.
For these larger problems, we set $\tolinner=5\cdot 10^{-9}$, $\tolouter=\tolIRKA=10^{-8}$.
The performance data of all methods are collected in Table~\ref{tbl:MIMO_large_problems_CPU}.
CPU times are in favor of our new approaches, with times that are from three to ten times lower.
These values are related to the lower number of outer iterations; depending on the coefficient
matrix sparsity, the number of \IRKA\ iterations is strictly correlated with the cost of linear system solves.
For instance, for {\sc rail20209}, solving linear systems in \IRKA\ takes 98\% of the total time, whereas
only 79\% for \RIRKA.

The required number of iterations for the truncated variant
\TRIRKA\ is higher than that for \RIRKA\, as expected, though memory requirements are significantly lower.
Indeed, we recall that \TRIRKA\ only works with $3\cdot 2r$ long vectors, as opposed to \RIRKA, which at the
end requires $(\#its+2) \cdot 2r$ vectors.

As for the previous example, Figure~\ref{figure_MIMO_rail20209} shows the convergence history (left) and
the location of the shifts at convergence (right), for the \textsc{rail20209} problem.

In Figure~\ref{figure_MIMO_rail20209_T}, again for the \textsc{rail20209}  data, we report the CPU time of each inner iteration,
together with the number of inner iterations for \RIRKA\ (left) and \TRIRKA\ (right). The maximum number of iterations (300)
is shown if the requested accuracy is not met. Nonetheless, the computed shifts are sufficiently good as shifts
for the next iteration, so that convergence is eventually achieved. As expected,
in the final stage of the outer process, the inner iteration number is minimal, as the  shifts gradually  stabilize.
This is a welcome event for \RIRKA, as the inner problem has growing dimensions, and thus becomes more expensive as
the outer iterations proceed. The different costs are emphasized on the y-axes.
\flushright $\diamond$
}
\end{example}
\vskip 0.1in

\begin{table}[htbp]
{\small
\centering
\caption{Elliptic operators and discretization size.\label{tab:matrices_L}}
\begin{tabular}{|l|r|l|}
		\hline
		Name& Size &Origin \\
		\hline
		\textsc{L10000} &10 000&$\mathcal {L}(u)= (\exp(-10xy)u_x)_x +(\exp(10xy)u_y)_y -(10(x + y)u)_x$  (\!\!\cite{MR2858340})\\ 
		\textsc{L10648} &10 648& $\mathcal {L}(u)= u_{xx}+u_{yy}+u_{zz} -10xu_x-1000yu_y-10u_z$  (\!\!\cite{MR2318706})\\
		\textsc{L160000} &\! 160 000&$\mathcal {L}(u) = {\rm div}(\exp(3xy)  {\rm grad} u) - 1/(x+y) u_x$ (\!\!\cite{MR2663662}) \\
		\hline
\end{tabular}
}
\end{table}

\begin{table}[htb]
\caption{Example~\ref{ex:largescale} Order $r=11$ \H2-optimal MOR of MIMO large scale problems.
For all methods, Option 1 was used as starting guess. \label{tbl:MIMO_large_problems_CPU}}
\centering
	{\scriptsize
\begin{tabular}{|c|crrrrrrr|}
			\hline
			&		&	\textsc{L10000}	&	\textsc{L10648}	&	\textsc{flow\_v0}	&	\textsc{flow\_v0.5}		&	\textsc{rail5177}	&	\textsc{rail20209}
			&	\textsc{t2dah} \\
			\hline
			\multirow{5}{*}{Info.}&	Size	&${\rm	10000	}$&${\rm	10648	}$&${\rm	9669	}$&${\rm	9669	}$&${\rm	5177	}$&${\rm	20209	}$&${\rm	11445	}$\\
			&	SYM	&${\rm	No	}$&${\rm	No	}$&${\rm	Yes	}$&${\rm	No	}$&${\rm	Yes	}$&${\rm	Yes	}$&${\rm	Yes	}$\\
			&	$B$	&${\rm	2}$&${\rm 2	}$&${\rm	1 }$&${\rm	1	}$&${\rm 7	}$&${\rm 7}$&${\rm 1}$\\
			&	$C$	&${\rm	2	}$&${\rm 2 }$&${\rm 5	}$&${\rm	5	}$&${\rm 6	}$&${\rm	6	}$&${\rm	7	}$\\
			\hline
			&	Method		&		&		&	&		&	&  &\\
			\hline
			
			\multirow{3}{*}{\#its}
&	\IRKA	&${\rm	34 	}$&${\rm	71 	}$&${\rm	26 	}$&${\rm	50 	}$&${\rm	131 	}$&${\rm	129 	}$&${\rm	24 	}$\\
&	\RIRKA	&${\rm	9 	}$&${\rm	8 	}$&${\rm	7 	}$&${\rm	6 	}$&${\rm	9 	}$&${\rm	9 	}$&${\rm	8 	}$\\
&	\TRIRKA	&${\rm	8 	}$&${\rm	8 	}$&${\rm	8 	}$&${\rm	6 	}$&${\rm	12 	}$&${\rm	13 	}$&${\rm	9 	}$\\

			\hline
			
			\multirow{3}{*}{$\linsolvers$}
&	\IRKA	&${\rm	748 	}$&${\rm	1562 	}$&${\rm	572 	}$&${\rm	1100 	}$&${\rm	2882 	}$&${\rm	2838 	}$&${\rm	528 	}$\\
&	\RIRKA	&${\rm	198 	}$&${\rm	176 	}$&${\rm	154 	}$&${\rm	132 	}$&${\rm	198 	}$&${\rm	198 	}$&${\rm	176 	}$\\
&	\TRIRKA	&${\rm	176 	}$&${\rm	176 	}$&${\rm	176 	}$&${\rm	132 	}$&${\rm	264 	}$&${\rm	286 	}$&${\rm	198 	}$\\

			\hline
			\multirow{3}{*}{CPU}
&	\IRKA	&${\rm	14.78 	}$&${\rm	143.36 	}$&${\rm	7.36 	}$&${\rm	27.39 	}$&${\rm	11.32 	}$&${\rm	63.05	}$&${\rm	9.20 	}$\\
&	\RIRKA	&${\rm	4.92 	}$&${\rm	19.13 	}$&${\rm	2.33 	}$&${\rm	3.88 	}$&${\rm	1.48 	}$&${\rm	5.41	}$&${\rm	3.71 	}$\\
&	\TRIRKA	&${\rm	4.18 	}$&${\rm	17.89 	}$&${\rm	2.55 	}$&${\rm	3.79 	}$&${\rm	1.93 	}$&${\rm	7.55	}$&${\rm	4.19 	}$\\

			\hline
			\multirow{3}{*}{$\finalsubspace$}	
&	\IRKA	&${\rm	11	}$&${\rm	11	}$&${\rm	11	}$&${\rm	11	}$&${\rm	11 	}$&${\rm	11 	}$&${\rm	11 	}$\\
&	\RIRKA	&${\rm	121 	}$&${\rm	110 	}$&${\rm	99 	}$&${\rm	88 	}$&${\rm	121 	}$&${\rm	121 	}$&${\rm	110 	}$\\
&	\TRIRKA	&${\rm	33 	}$&${\rm	33 	}$&${\rm	33 	}$&${\rm	33 	}$&${\rm	33 	}$&${\rm	33 	}$&${\rm	33 	}$\\

	\hline
			
\end{tabular}
	}
\end{table}


In all previous experiments, Option 1 was used to compute the initial bases $V$ and $W$. 
This entails approximating the $2r$ smallest eigenvalues of the given matrices. For certain large
problems, this costs may be excessive, even for very lose stopping tolerances. 
In the following example we explore the possibility of starting with random vectors,
corresponding to Option 3 in the algorithm. 

\begin{figure}[htb]
 	\caption{Example~\ref{ex:largescale_rand}  Convergence behaviors of the algorithms for \textsc{GasSensor} and \textsc{T3DH}, MIMO, \texttt{ $r=16$}.
Left: GasSensor. Right: T3DH.
 			\label{figure:gas_and_t3dh} }
 			\includegraphics[width=2.5in]{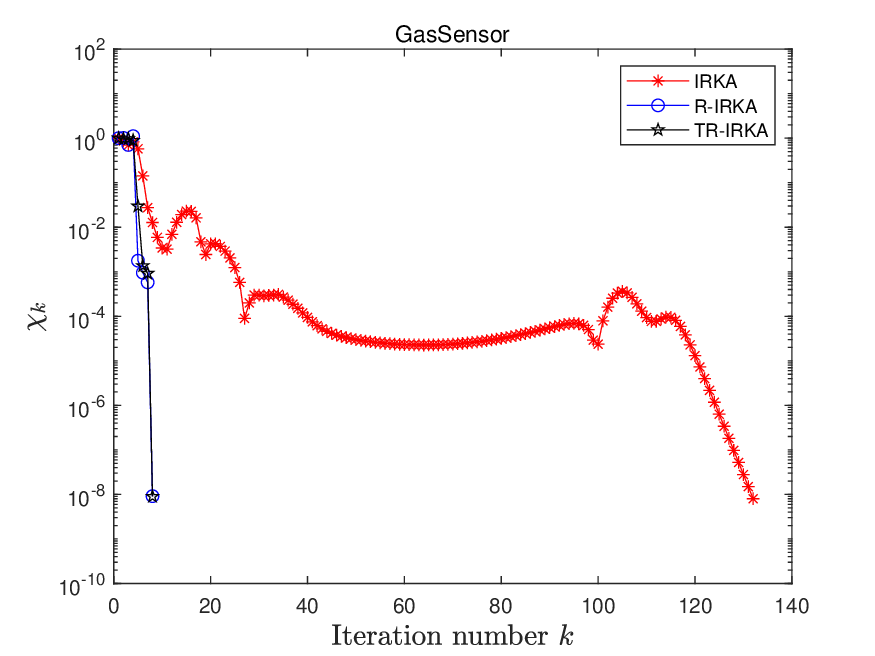} 
 			\includegraphics[width=2.5in]{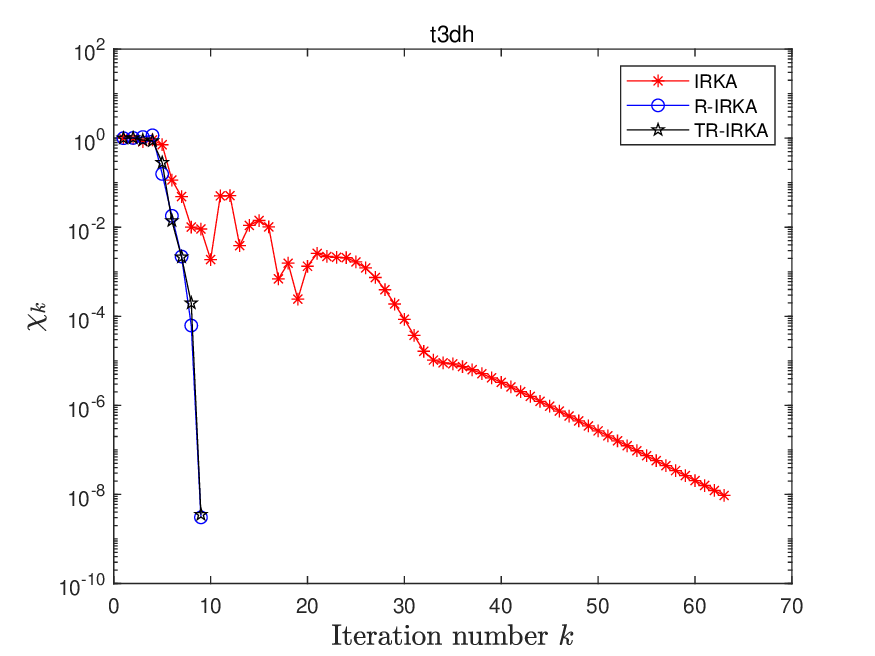} 
 \end{figure}

\vskip 0.1in
\begin{example}\label{ex:largescale_rand}
{\rm 
We consider a dataset of larger problems from the Oberwolfach collection, together with
the data in \cite[Example 4.4]{MR4082008}, where the coefficient matrix has Toeplitz structure.
This latter one appeared to be quite a challenging eigenproblem.  
As starting guess, we consider $V, W$ obtained with Option 3 in the corresponding algorithms.
The performance results are presented in Table~\ref{tbl:MIMO_example_2}, where the same
stopping tolerances as in the previous example have been used.
For \textsc{Toeplitz} and \textsc{L160000} (from Table~\ref{tab:matrices_L}), we set $E=I$ and obtain 
$B$ and $C$ by \texttt{randn(`state',0);} \texttt{B=randn(n,2);} \texttt{randn(`state',10);} \texttt{C=randn(n,2)}.

The CPU times in Table~\ref{tbl:MIMO_example_2} confirm the competitiveness of the
new algorithms over \IRKA\ when solving large-scale problems.
For completeness, the convergence history for
\textsc{GasSensor} and \textsc{T3D} is displayed in Figure~\ref{figure:gas_and_t3dh}.
For all methods the large majority of the computational efforts 
(up to about \%99) is focused on solving linear equations
with long vectors, 
we can observe a significant reduction in overall CPU time whenever the number of outer iterations 
is largely reduced from \IRKA\ to \RIRKA.
The \textsc{Toeplitz} problem is an exception: thanks to the structure, solving with the Toeplitz matrix is
not very expensive, in spite of the dimension. In this case, all other costs (orthogonalization, reduced problem, etc.)
become more relevant. This explains CPU times that are not as different for the considered methods as
the number of outer iterations would predict. 

Finally,  we comment on the convergence rate.
In the majority of the observed examples \IRKA\ converges linearly \cite{MR2924212,10740341},
though the iteration numbers can vary considerably for different data.
For instance, on \textsc{Filter3D}, \IRKA\ only requires  $30$ iterations.
On the other hand,
 on \textsc{L160000}, \IRKA\ requires  $95$ iterations. In all cases, \IRKA\ shows a transient period during
which the next generated shift set does not improve the approximation to the optimal set, so that an almost
stagnating convergence curve can be observed;
see, e.g., the plots in Figure~\ref{figure:gas_and_t3dh}.
In \RIRKA, this transient phase is mostly taken care of by the reduced problem, so that stagnation
of the outer iteration is minimal, and superlinear final convergence can be observed.
\flushright $\diamond$
}
\end{example}
\vskip 0.1in

\begin{table}[h]
\centering
\caption{Example~\ref{ex:largescale_rand}. Order $r=16$ \H2-optimal MOR of large scale  problems(MIMO).
Starting bases obtained with Option 3. \label{tbl:MIMO_example_2}}
	\footnotesize 
\begin{tabular}{|c|rrrrrrrr|}
\hline
		&	&	{\scriptsize\textsc{Toeplitz}}	&{\scriptsize\textsc{T3DL}}&	{\scriptsize\textsc{L160000}}			&{\scriptsize\textsc{rail79841}}	&	{\scriptsize\textsc{filter3D}}	&	{\scriptsize\textsc{GasSensor}} &	{\scriptsize\textsc{T3DH}}\\
		
\hline
\multirow{4}{*}{Info.}
		&	Size	&${\rm	200000	}$&${\rm	20360	}$&${\rm	160000	}$&${\rm	79841	}$&${\rm	106437	}$&${\rm	66917 	}$&${\rm	79171 	}$\\
		&	SYM	&${\rm	No	}$&${\rm	Yes	}$&${\rm	No	}$&${\rm	Yes	}$&${\rm	Yes	}$&${\rm	Yes	}$&${\rm	Yes	}$\\
		&	$B$	&${\rm	2	}$&${\rm	1	}$&${\rm	2	}$&${\rm	7	}$&${\rm	1	}$&${\rm	1 	}$&${\rm	1 	}$\\
		&	$C$	&${\rm	2	}$&${\rm	7	}$&${\rm	2	}$&${\rm	6	}$&${\rm	5	}$&${\rm	28 	}$&${\rm	7 	}$\\
		
\hline
			&	Method		&		&		&	&		&	&  & \\
\hline
		
\multirow{3}{*}{\#its}
&	\IRKA	&${\rm	18 	}$&${\rm	37 	}$&${\rm	95 	}$&${\rm	104 	}$&${\rm	30 	}$&${\rm	132 	}$&${\rm	63 	}$\\
&	\RIRKA	&${\rm	5 	}$&${\rm	7 	}$&${\rm	8 	}$&${\rm	11 	}$&${\rm	8 	}$&${\rm	8 	}$&${\rm	9 	}$\\
&	\TRIRKA	&${\rm	5 	}$&${\rm	8 	}$&${\rm	8 	}$&${\rm	16 	}$&${\rm	8 	}$&${\rm	8 	}$&${\rm	9 	}$\\

		\hline
		
		\multirow{3}{*}{$\linsolvers$}
&	\IRKA	&${\rm	576 	}$&${\rm	1184 	}$&${\rm	3040 	}$&${\rm	3328 	}$&${\rm	960 	}$&${\rm	4224 	}$&${\rm	2016 	}$\\
&	\RIRKA	&${\rm	160 	}$&${\rm	224 	}$&${\rm	256 	}$&${\rm	352 	}$&${\rm	256 	}$&${\rm	256 	}$&${\rm	288 	}$\\
&	\TRIRKA	&${\rm	160 	}$&${\rm	256 	}$&${\rm	256 	}$&${\rm	512 	}$&${\rm	256 	}$&${\rm	256 	}$&${\rm	288 	}$\\

		\hline
		\multirow{3}{*}{CPU}
&	\IRKA	&${\rm	16.1 	}$&${\rm	514.4 	}$&${\rm	1622.2 	}$&${\rm	310.3 	}$&${\rm	1135.8 	}$&${\rm	5178.5 	}$&${\rm	7667.2 	}$\\
&	\RIRKA	&${\rm	7.5 	}$&${\rm	95.5 	}$&${\rm	144.6 	}$&${\rm	49.1 	}$&${\rm	319.5 	}$&${\rm	449.8 	}$&${\rm	1581.2 	}$\\
&	\TRIRKA	&${\rm	6.8 	}$&${\rm	104.8 	}$&${\rm	141.8 	}$&${\rm	58.1 	}$&${\rm	291.4 	}$&${\rm	454.9 	}$&${\rm	1430.2 	}$\\

		\hline
		\multirow{3}{*}{$\finalsubspace$}
		&	\IRKA	&${\rm	16	}$&${\rm	16	}$&${\rm	16	}$&${\rm	16	}$&${\rm	16	}$&${\rm	16	}$&${\rm	16	}$\\		
&	\RIRKA	&${\rm	112 	}$&${\rm	144 	}$&${\rm	160 	}$&${\rm	208 	}$&${\rm	160 	}$&${\rm	160 	}$&${\rm	176 	}$\\
&	\TRIRKA	&${\rm	48 	}$&${\rm	48 	}$&${\rm	48 	}$&${\rm	48 	}$&${\rm	48 	}$&${\rm	48 	}$&${\rm	48 	}$\\
\hline
\end{tabular}
\end{table}

\section{Conclusions}\label{sect:conclusion}
We have developed a novel rational Krylov subspace method for approximating
the optimal shifts giving \H2-optimal model order reduction.
For large scale problems, our strategy \RIRKA\ accumulates the rational spaces
associated with the most recently computed sets of shifts, so that the classical
\IRKA\ can exploit a richer space.

Our computational results are promising. Indeed,
our experiments with various benchmark problems showed that
the number of \RIRKA\ iterations is significantly lower than that of \IRKA.
For the large-scale problems, this difference corresponds to a dramatic
CPU time reduction, for similar quality results in terms of transfer function approximation
and optimal parameter selection.

Understanding the theory supporting these results requires a deeper analysis
of the procedure, which will require a separate project.

\section*{Appendix 1}
In this appendix we report the algorithm of
 the classical IRKA \cite{MR2421462} for the 
\H2-optimal MOR in Algorithm~\ref{alg:IRKA}.
Matlab notation is used whenever possible.

\begin{algorithm}[tbh]
	\caption{An Iterative Rational Krylov Algorithm (IRKA) for the  \H2-optimal MOR \cite[Algorithm 4.1]{MR2421462}.
		\label{alg:IRKA}}

	\begin{algorithmic}[1]
		\REQUIRE $A,E\in \mathbb{C}^{n \times n},b,c\in \mathbb{C}^{n \times 1}$.  
		Order $r$.	 Option number at Line 1.
		
		\STATE Get initial bases of the projection subspaces:	
		\begin{eqnarray*}
				{\rm Option \, 1:} &&
					[V, \sim ]=\texttt{eigs}(A,E,r,0);
				V = \orth(V), W = V\\
				{\rm Option \, 2:} &&
				{\rm Given}: \mathbb{S}_0, \mbox{Create } 
{V, W 
\mbox{ orth. basis of } \RK_r({A},{E}, b,\mathbb{S}_0), \RK_r({A}^H,{E}^H, c,\overline{\mathbb{S}}_0)}\\
				%
				{\rm Option \, 3:}&& 
				V = \orth(\texttt{randn}(n,r)); W = \orth(\texttt{randn}(n,r)).
		\end{eqnarray*}
		
		\STATE 
		Set  $\mathbb{S}_0=\texttt{ones}(r,1)$ in Option $1$ or Option $3$.
		Set $\mathbb{S}_0=\texttt{sort}(\mathbb{S}_0)$ in  Option $2$.
		
		\FOR{$k=1,2, \ldots, k^{\max} $}
	\STATE $\Lambda=\texttt{eig}((W^HEV)^{-1}W^HAV).$
		\STATE $\mathbb{S}_k=-\texttt{conj}(\Lambda)$, $\mathbb{S}_k=\texttt{sort}(\mathbb{S}_k ).$
		\STATE  $V=\orth(\RK({A},{E}, b,\mathbb{S}_k,r)), W=\orth(\RK({A}^H,{E}^H, c,\overline{\mathbb{S}_k},r)).$
		\STATE $\chi_k^{\rm rel}:=
		\|\mathbb{S}_k-\mathbb{S}_{k-1}\|_2 /\|\mathbb{S}_k\|_2.$
		\STATE  {\bf if}
		$\chi_k^{\rm rel}< \texttt{tol}_3$  {\bf then} {\bf break};  {\bf end if}
		
		\ENDFOR
		
	   \ENSURE   \H2-optimal MOR  shifts and subspaces: $\mathbb{S}_k$, $\Range({V})$, $\Range({W})$.
	\end{algorithmic}
\end{algorithm}


\section*{Acknowledgments}
The authors thank Serkan Gugercin and Zlatko Drma\v{c} for making 
their Matlab codes of \IRKA\
available, and Chris Beattie for an insightful conversation on \IRKA.
The work of Y.L. is supported by the National Natural Science 
Foundation of China(NSFC-12101508).
V.S. is a member of the INdAM Research
Group GNCS.  Its continuous support is fully acknowledged.

\bibliography{math20250115}

\end{document}